\documentclass[12pt,reqno]{amsart}
\usepackage{amssymb}
\usepackage{graphicx}
\input xy
\xyoption{all}
\title[A Zero-One law]{A Zero-One Law for Random Subgroups of some Totally Disconnected Groups.}
\author{Yair Glasner}
\email{yairgl@math.bgu.ac.il}
\address{Department of Mathematics, Ben Gurion university of the Negev, Be'er Sheva 84105, ISRAEL.}

\newtheorem*{theorem}{Theorem}

\newtheorem*{proposition}{Proposition}
\newtheorem*{corollary}{Corollary}
 \theoremstyle{definition}
\newtheorem*{question}{Question}
\newtheorem*{definition}{Definition}
\newtheorem*{remark}{Remark}

\newtheorem*{examples}{Examples}

\newcommand{\Z}{{\mathbf{Z}}}
\newcommand{\R}{{\mathbf{R}}}
\newcommand{\C}{{\mathbf{C}}}
\newcommand{\Q}{{\mathbf{Q}}}
\newcommand{\F}{{\mathbf{F}}}

\renewcommand{\a}{{\mathbf{a}}}

\newcommand{\Oc}{{\mathcal{O}}}

\newcommand{\arrow}{\rightarrow}
\newcommand{\Hom}{{\operatorname{Hom}}}

\newcommand{\chr}{{\operatorname{char}}}

\newcommand{\Aut}{{\operatorname{Aut}}}
\newcommand{\Comm}{{\operatorname{Comm}}}

\newcommand{\tr}{{\operatorname{\mathcal{TR}}}}

\newcommand{\PSL}{{\operatorname{PSL}}}
\newcommand{\PGL}{{\operatorname{PGL}}}
\newcommand{\SL}{{\operatorname{SL}}}

\newcommand{\GL}{{\operatorname{GL}}}
\newcommand{\SU}{{\operatorname{SU}}}
\newcommand{\sll}{{\operatorname{\mathfrak{sl}}}}

\newcommand{\defeq}{\stackrel{\operatorname{def}}{=}}

\newcommand{\Ell}{{\operatorname{Ell}}}
\newcommand{\Hyp}{{\operatorname{Hyp}}}

\newcommand{\Ad}{{\operatorname{Ad}}}

\newcommand{\compos}{\circ}

\begin{document}
\bibliographystyle{alpha}

\begin{abstract}
Let $A$ be a locally compact group topologically generated by $d$ elements and let $k > d$. Consider the action, by pre-composition, of $\Gamma = \Aut(F_k)$ on the set of marked, $k$-generated, dense subgroups $D_{k,A} := \{ \eta \in \Hom(F_k,A) \ | \ \overline{\langle \phi(F_k) \rangle} = A\}.$ We prove the ergodicity of this action for the following two families of simple, totally disconnected locally compact groups: 
\begin{itemize}
\item $A = \PSL_2(K)$ where $K$ is a non-Archimedean local field (of characteristic $\ne 2$),
\item $A = \Aut^{0}(T_{q+1})$ - the group of orientation preserving automorphisms of a $q+1$ regular tree, for $q \ge 2$.
\end{itemize} 
In contrast, a recent result of Minsky's shows that the same action is not ergodic when $A = \PSL_2(\R)$ or $A = \PSL_2(\C)$. Therefore if $K$ is a local field (with $\chr(K) \ne 2$) the action of $Aut(F_k)$ on $D_{k,\PSL_2(K)}$ is ergodic, for every $k > 2$, if and only if $K$ is non-Archimedean. 

Ergodicity implies that every ``measurable property'' either holds or fails to hold for almost every $k$-generated dense subgroup of $A$.
\end{abstract}

\maketitle
\section{Introduction}
\subsection{Invariants of $k$-generated subgroups}
Throughout the paper $A$ will denote a locally compact group, topologically generated by $d = d(A)$ elements, with Haar measure $\mu$. We will fix an integer $k > d$ and consider the collection of $k$-tuples
$A^k \cong \Hom(F_k,A)$, this group comes with the (possibly infinite) product Haar measure $\mu^k$. 

Since every locally compact group automatically combines the structure of a group with that of a measure space, the following question arises in a very natural way:
\begin{question} \label{q:Wie} \nonumber
What are the properties of a subgroup of $A$ generated by a generic $k$-tuple of elements?
\end{question}
By generic we mean excluding a subset of $\mu^k$-measure zero. One obvious invariant is the closure $\overline{\eta(F_k)} < A$ of such a subgroup. In fact, because of our freedom to replace $A$ by a closed subgroup, it is natural to restrict our attention to the set $$D = D_{k,A} \defeq \{\eta \in \Hom(F_k,A) \ | \ \overline{\eta(F_k)} = A\},$$ of $k$-tuples generating a dense subgroup of $A$.

Given any open set $U \subset A$ and a specific $w \in F_k$ the collection of homomorphisms $\{ \phi \in \Hom(F_k,A) \ | \ \phi(w) \in U \}$ is open. The set $D_{k,A}$ is therefore $G_{\delta}$ and hence Borel measurable, whenever $A$ is second countable. It might happen though that $D_{k,A}$ is a nullset but for many nice groups $A$ this is not the case and often $D_{k,A}$ is even open. This is the case, for example, for all connected Lie groups as well as for the two families of groups, which are at the focus of this paper, mentioned in the abstract. 

\subsection{The Goldman-Wiegold question}
Let $\Gamma = \Aut(F_k)$ be the automorphism group of the free group on $k$ generators. $\Gamma$ acts on $\Hom(F_k,A)$ by pre-composition. When $A$ is unimodular this action is measure preserving. In general, this action need not be measure preserving but it always preserves the measure class of Haar measure and in particular takes nullsets to nullsets. Clearly $\eta(F_k) = \eta \compos \gamma (F_k)$ for every $\gamma \in \Gamma$ and in particular $D$ is invariant under the $\Gamma$ action. We will be interested in the following
\begin{question} (The Goldman-Wiegold question)  \label{Q:GW}
For which groups $A$ is the action of $\Gamma = \Aut(F_k)$ ergodic on $D = D_{k,A}$ for every $k > d(A)$?
\end{question}
When this action is ergodic almost all $k$-generated dense subgroups of $A$ look alike: every measurable property either holds or fails to hold for almost all of them together. A group property $P$ is {\it{measurable}} if the set $\{\eta \in D_{k,A} \ | \ \eta(F_k) {\text{ has the property P}}\}$ is measurable. One example of such a property is the isomorphism type of the group. Other properties, such as the existence of a spectral gap treated in \cite{Fisher:spectral_gap} in this setting, are more context specific. 

\subsection{The main theorem} Let $A = \PSL_2(K)$ be the projective special linear group over a non-Archimedean local field $K$ of characteristic $p \ne 2$; or $A = \Aut^0(T)$ the group of orientation preserving automorphisms of a regular tree. Both classes of groups are totally disconnected and topologically generated by two elements. The latter fact follows as a special case of Theorem \ref{thm:SL_gen} below for $A = \PSL_2(K)$ and Theorem \ref{thm:AG} for $A = \Aut^{0}(T)$. Furthermore both groups are simple and in particular unimodular. We give an affirmative answer to the Goldman-Wiegold question for these groups:
\begin{theorem} \label{thm:GW} (Main theorem)
Let $A$ be one of the groups described above. Then for every $k \ge 3$, the set $D_{k,A}$ is of infinite Haar measure in $\Hom(F_k,A)$ and the  action of $\Gamma = \Aut(F_k)$ on $D_{k,A}$ is ergodic.
\end{theorem} 
Excluding certain easy cases such as finitely generated Abelian groups, these are the first examples outside the world of compact groups for a positive answer to the Goldman-Wiegold question.

For any measurable group property $P$ it is now interesting to understand whether $P$ holds or fails to hold for almost every $k$-generated subgroup of $A$.  Two such properties are treated in \cite{AG:boundary}, where it is proved that almost every $k$-generated subgroup of $\Aut(T)$ is free and acts essentially freely on the boundary of the tree. Here an action of a group on a measure space is called {\it{essentially free}} if the set of fixed points of every non-trivial element is a nullset. 

\subsection{Minsky's results in the Archimedean case}
While this paper was being refereed, Yair Minsky has been addressing the same ergodicity question, for $\SL_2(K)$ over the Archimedean local fields $K=\R$ and $K=\C$. To my great surprise the outcome of these investigations shows that the action of $\Aut(F_n)$ is  {\it{not ergodic}} in the Archimedean case:
\begin{theorem} (Yair Minsky)
Let $A = \PSL_2(\R)$ or $A = \PSL_2(\C)$ then for any $k \ge 3$ then there are two disjoint $\Aut(F_k)$-invariant subsets with non-empty interior $$PS, R \subset D_{k,A}$$
\end{theorem}
The set $PS$ is also referred to as the {\it{primitive stable}} part and $R$ as the redundant part. Minsky shows, in addition, that the action of $\Aut(F_n)$ on $PS$ is properly discontinuous. It was recently proved by Avni, Gelander and Minsky 
\cite{AGM:Redundant} that the action on the redundant part is ergodic. Conjecturally $PS \cup R$ is conull in $D_{k,A}$. 

Combining the results of Minsky with these of the current paper gives rise to the following:
\begin{corollary}
Let $K$ be a local field with $\chr(K) \ne 2$ and $k \ge 3$. Then the action of $\Aut(F_k)$ on $D_{k,A} := \left \{\phi \in \Hom(F_k,A) \ | \ \overline{\phi(F_k)} = A \right\}$ is ergodic if and only if $K$ is non Archemedian. 
\end{corollary}
\subsection{History of the Goldman-Wiegold question.} We will describe the background and history of the Goldman-Wiegold question only briefly, referring the readers to an upcoming comprehensive survey by Alex Lubotzky \cite{Lub:Aut_Fn} dedicated to this subject.

For finite simple groups Question \ref{Q:GW} is a well known conjecture, attributed to J. Wiegold. By the classification of finite simple groups every such group is $2$-generated so the conjecture in this case asserts the transitivity of the action of $\Aut(F_k)$ on $k$-tuples of generators. Gilman \cite{Gilman:SL2p} proved Wiegold's conjecture for $A = \SL_2(F_p)$ when $p$ is prime. Evans \cite{Evans:Wiegold_2} established Wiegold's conjecture for $A = \SL_2(F_{2^n})$ and for $A = \operatorname{Sz} (2^{2m-1})$. Allowing for slightly larger values of $k$, Garion \cite{Garion:PSL2} proved transitivity of the $\Aut(F_4)$ action on quadruples of generators for $\SL_2(F_q)$ where $F_q$ is an arbitrary finite field. For general finite simple groups of Lie type, Avni and Garion \cite{AG:bounded_rk} prove the transitivity of the $\Aut(F_k)$ action on $k$-tuples of generators, where $k = k(r)$ depends only on the Lie rank of the group. For finite solvable groups an affirmative answer to Question \ref{Q:GW} is given by Dunwoody in \cite{Dunwoody:Nielsen}. 

Independently of all the activity in the finite case, W. Goldman \cite{Goldman:Ergodic} conjectured an affirmative answer to Question \ref{Q:GW} for compact connected Lie groups and proved the special case $A = \SU_2$. For compact connected groups it turns out that almost all $k$-tuples generate a dense subgroup, so up to nullsets $D_{k,A} = \Hom(F_k,A)$. Goldman's conjecture was proved in full generality in a beautiful paper by T. Gelander \cite{Gelander:deff} establishing the ergodicity of $\Aut(F_k)$ on $\Hom(F_k,A)$ for every compact connected Lie group $A$ and any $k \ge 3$. 

As observed above, ergodicity of $Aut(F_k)$ in the above examples, implies a zero-one law for subgroups. But usually it is not easy to find measurable group properties for which the zero-one statement is interesting and unknown. One successful application of Goldman's theorem is due to David Fisher \cite{Fisher:spectral_gap} who observed that $k$ random elements $a_1,a_2,\ldots,a_k \in \SU_2$ almost surely either admit, or fail to admit a spectral gap, in their action on 2-sphere. An open conjecture, known as the spectral gap conjecture, due to Lubotzky Phillips and Sarnak \cite{LPS:SphereI,LPS:SphereII}, asserts that that a spectral gap is indeed generic. 

Prior to the recent results of Minsky, the only source of examples for situations where ergodicity does fail came from the world of infinite discrete groups. But in that setting many examples were found - even for metabelian groups, see for example \cite{Evans:more_gen,Evans:more_genII}.  Ergodicity also fails miserably once we consider the case $k = d$. See \cite{GS:T_sys} for the most recent results on how transitivity fails for finite simple groups when $k = 2$. In the case that interests us - the group $\PSL_2(K)$ - the trace of the commutator, $\operatorname{tr}([a,b])$ is a measurable function on $A^2$ invariant under the $\Gamma$ action, which immediately contradicts ergodicity. 

It was mentioned to me by Alex Lubotzky that Question \ref{Q:GW} probably assumes rather different nature in the cases where $D_{k,A}$ is a nullset for arbitrarily large values of $k$. By a theorem of Kantor and Lubotzky \cite{KL:prob_gen} this is indeed the case when $A = \hat{F_2}$ is the free profinite group: the probability that $k$ random elements of $A$ generate an open subgroup is zero for every $k$. 

We conclude this section by mentioning that there are many variations on Question \ref{Q:GW}. For example some authors consider the diagonal action of $\Aut(F_k) \times \Aut(A)$ on $\Hom(F_k,A)$, the orbits of this action are referred to as as T-systems in the case where $A$ is finite. 

\subsection*{Thanks} This paper was originally written with a different, more geometric, proof for Theorem \ref{thm:SL_gen} that was valid only in characteristic zero (i.e. when $K$ is a finite extension of $\Q_p$). I would like to thank Alex Lubotzky who suggested the method for generalizing these results to positive characteristic. I would also like to thank Tsachik Gelander for many interesting discussions on this subject. Finally, many thanks to the two referees who made many helpful remarks on the previous version of this paper. This research was partially funded by ISF grant 888/07 and by BSF grant 2006222. 

\section{Notation and Preliminary results}
\subsection{Nielsen Transformations}
It is a famous theorem of Nielsen that the group $\Gamma = \Aut(F_k)$ is finitely generated. The following, particularly nice set of generators, is referred to as Nielsen transformations:
\begin{itemize}
\item For every $1 \le i \ne j \le k$ we define:
\begin{eqnarray*}
R^{\pm}_{i,j}(x_1,\ldots,x_k) & = & (x_1,\ldots,x_{j-1}, x_jx_i^{\pm}, x_{j+1}, \ldots x_k) \\
L^{\pm}_{i,j}(x_1,\ldots,x_k) & = & (x_1, \ldots,x_{j-1}, x_i^{\pm}x_j, x_{j+1}, \ldots x_k) 
\end{eqnarray*}
\item For every $1 \le i < j \le k$ define:
\begin{eqnarray*}
T_{i,j} (x_1, \ldots, x_i, \ldots, x_j, \ldots x_k) & = & (x_1, \ldots, x_j, \ldots, x_i, \ldots, x_k)
\end{eqnarray*}
\end{itemize}

The Schreier graph coming from the action of $\Gamma$ on the set of $k$-tuples of generators $D \subset A^k$ with respect to this specific choice of generating set is known as the {\it product replacement graph}. In the case where $A$ is finite and simple the  Wiegold conjecture asserts the connectedness of the product replacement graph. We refer the reader to \cite{LP:PRA} and the references therein for more on this point of view. 

\subsection{Classification of tree automorphisms.} Both groups that play central role in this paper $\Aut^{0}(T)$ and $\PSL_2(K)$ come with a natural action on a tree. For the former the action is given explicitly while the latter acts on its associated Bruhat-Tits tree. 

We recall here the classification of tree automorphisms according to their dynamical properties, for complete proofs see (\cite{Serre:Arbres}).
Let $T = T_{q+1}$ be a $(q+1)$-regular tree $\Aut(T)$ its automorphism group and $\Aut^0(T)$ the index two subgroup that preserves the natural $2$-coloring of the vertices. For $a \in \Aut(T)$, we denote by $\ell(a) = \min \{d_T(x,ax) \ | \ x \in T \}$ and by $X(a) \defeq \{x \in T \ | \ d_T(x,ax) = \ell(a) \}$, where the minimum is taken over all points in the {\it{geometric realization of the tree}}. If $\ell(a) = 0$ and $X(a)$ contains a vertex then $a$ is called {\it{an elliptic element}}. In this case $X(a)$ is the subtree of points fixed by $a$. If $\ell(a) = 0$ but $X(a)$ contains no vertices then $a$ is called {\it{ an inversion}} and $X(a)$ is a unique point in the center of one of the geometric edges. If $\ell(a) > 0$ then $a$ is known as a {\it{hyperbolic element}}. In this case $X(a)$ is an infinite geodesic on which $a$ acts as a translation of length $\ell(a)$. 
\begin{definition}
We denote by $\Ell$ the set of elliptic elements of $\Aut(T)$ and by $\Hyp$ the set of hyperbolic elements. Inside $\PSL_2(K)$ the same symbols $\Ell$ (resp. $\Hyp$) will denote elements that are elliptic (resp. hyperbolic) in their action on the corresponding Bruhat-Tits tree. 
\end{definition}
We will not deal much with inversions, since our interest is focused on two groups that do not contain inversions: $\Aut^{0}(T)$ and $\PSL_2(K)$. 

\subsection{The Nielsen method.}
In a key point in the proof we will make use of the following theorem due to Richard Weidmann.
\begin{theorem} \cite[Theorem 7]{Weidmann:Nielsen} \label{thm:Weidmann}
\label{weidman} Let $A = \Aut(T_{q+1})$ be the automorphism group of a regular tree and $\eta \in \Hom(F_k,A)$. Then either $\eta(F_k)$ is a free group acting freely on $T$ or there is an element $\gamma \in \Gamma = \Aut(F_k)$ such that $\eta \compos \gamma (x_1)$ is elliptic.
\end{theorem}
The proof of the theorem uses the so called Nielsen method - which carefully follows the dynamical properties of automorphisms and how they transform under Nielsen transformations. This method was first used by Nielsen to show that a finitely generated subgroup of a free group is free - prior to the proof of the Nielsen-Schreier theorem in its complete generality.

\subsection{A subgroup of $\Aut(T)$ generated by an elliptic and a hyperbolic is generically dense.}
The generic properties of a group generated by $k$ elements in the group $A = \Aut^0(T)$, was studied extensively in \cite{AG:Generic,AG:boundary}. In particular we have shown in that paper that for $k = 2$ the set $\Ell \times \Hyp$ of elliptic-hyperbolic pairs is essentially contained in $D$ - the set of pairs that topologically generate $\Aut(T)$:
\begin{theorem} \cite{AG:Generic} \label{thm:AG} 
Let $T = T_{q+1}$ be the $q + 1$-regular tree for $q \ge 2$ and $A = \Aut^0(T)$ its automorphism group. Consider the set
$$\Hyp \times \Ell = \{(a,b) \in A^2 \ |  \ a {\text{ is hyperbolic and }} b {\text{ elliptic}} \}.$$  
Then $\Hyp \times \Ell$ is of infinite Haar measure and almost every $(a,b) \in \Hyp \times \Ell$ generates a dense subgroup of $A$.  
\end{theorem}
In Section \ref{sec:SL_gen} below we establish the analogous theorem  for the group $A = \PSL_2(K)$, where $K$ is a non-Archimedean local field. 

\subsection{Notation pertaining to local fields} \label{sec:notation}
$K$ will denote a non-Archimedean local field, i.e. either a finite extension of $\Q_p$ or isomorphic to $\F_q((t))$ for some prime power $q = p^l$. We denote by $\Oc$ the ring of integers in $K$, by $\pi$ a uniformizer and $k = \Oc/\pi \Oc = F_{p^l}$ the residue field. {\it To avoid technical complications we will always assume that $p = \chr(K) \ne 2$. }

Let $A = \PSL_2(K)$ and $T$ the corresponding Bruhat-Tits tree. $T$ is a $(q+1)$-regular tree on which $A$ acts by isomorphisms. We will refer to elements of $A$ as elliptic (resp. hyperbolic) if they are elliptic (resp. hyperbolic) in their action on $T$. Let $\tr: \PSL_2(K) \arrow K$ denote the trace of the adjoint 3-dimensional representation on the Lie algebra $\sll_2(K)$.

\subsection{Pink's classification of compact subgroups.}
Richard Pink in \cite{Pink:Compact},  gives a classification theorem of compact subgroups of $\GL_n(K)$ where $K$ is a local field. The main step of the proof, is understanding compact groups whose Zariski closure is an absolutely simple adjoint connected algebraic group. The simplest case and the one which will be of interest for us here, is the one of a compact Zariski dense subgroup in $\PSL_2(K)$. In this specific case, Pinks theorem reads:
\begin{theorem}(See \cite[Theorem 0.7]{Pink:Compact}) \label{thm:Pink}
Let $K$ be a non-Archimedean local field of characteristic $\ne 2$ and let $C < \PSL_2(K)$ be a Zariski dense, compact, subgroup. Then there exists a model of $\PSL_2$ over a closed subfield $E < K$, such that $C$ is an open subgroup of $H(E)$.
\end{theorem}
\begin{remark}
Note that Pink's theorem was well known for local fields of characteristic zero. The main novelty (and the most difficult case) in his proof is in positive characteristic fields.
\end{remark}
\begin{remark}
Pink theorem above works more generally for any absolutely simple, connected adjoint group $G$ under the additional assumption that the adjoint representation is irreducible. There is a more complete, but also more technical, form of Pink's theorem \cite[Theorem 0.2]{Pink:Compact} that does not require that the adjoint representation be irreducible.

In our setting the adjoint representation of $\PGL_2(K)$ is irreducible if and only if $\chr(K) \ne 2$. Indeed, since the adjoint representation is 3-dimensional, if it is reducible then it contains an invariant vector. Namely a matrix $A \in \sll_2(K)$ that commutes with all elements of $\PSL_2(K)$. It is easy to verify that such a matrix $A$ must be scalar; and a scalar matrix has trace zero if and only if the field has characteristic $2$. 
\end{remark}

\begin{examples} of non-open Zariski dense compact subgroups.
\begin{itemize}
\item {\it{Examples coming from subfields.}} $\PSL_2(\Oc) \cap \PSL_2(F) < \PSL_2(K)$, where $F<K$ is a local subfield. For example $\PSL_2(\F_p[[t]]) < 
PSL_2(\F_q((t)))$ with $q = p^l$. Or $\PSL_2(\F_p[[t^m]]) < \PSL_2(\F_p((t)))$ for some $m \in \Z$.   
\item {\it{Examples coming from compact forms.}}Assume that $Q$ is a quaternion algebra that ramifies over some subfield$F< K$ but splits over $K$. Let $H = Q^{*}/\mathcal{Z}Q^{*}$ be the projective group of invertible elements. Since $Q$ splits over $K$ we have an isomorphism $$\phi: H (K) \arrow \PSL_2(K)$$ and $\phi(H(F))$ is a compact Zariski dense subgroup of $\PSL_2(K)$. 
\end{itemize}
\end{examples}
Pink's theorem (often in its more general form \cite[Theorem 0.2]{Pink:Compact}) is extremely useful. In the following two sections we review two applications of this theorem, adapted to our specific cause.

\subsection{A theorem of Barnea-Larsen on random generation} It follows from Pink's Theorem \ref{thm:Pink} that if $C$ is a compact subgroup of $\PSL_2(K)$ such that
\begin{enumerate}
\item \label{itm:Zdense} $C$ is Zariski dense and
\item \label{itm:trace} $\tr(C)$ is not contained in any closed subfield of $K$
\end{enumerate}
then $C$ is open. In \cite[Theorem 3.1]{BL:Random_gen} Barnea and Larsen prove that if $U < A$ is compact open, then two Haar random elements $x,y \in U$ almost surely generate a group that satisfies both of the above conditions and is hence, almost surely open. In fact their proof gives the following stronger statement, which we will use later:
\begin{proposition} \label{prop:BL_conditions}
Let $K$ be a non-Archimedean local field of characteristic $\ne 2$, $A = \PSL_2(K)$ and $\mu$ Haar measure on $A$. Then
$$\mu \left\{ (a,b) \in A^2 \ | \ \langle x,y \rangle {\text{ does not satisfy (\ref{itm:Zdense}) or (\ref{itm:trace}) above }} \right\} = 0$$
\end{proposition}
\begin{proof}
The fact that the collection $$\left\{(a,b) \in A^2 \ | \ \langle a, b \rangle {\text{ is not Zariski dense}} \right\}$$ is of measure zero, and in fact contained in a countable union of smaller dimension varieties, is proved in \cite[Lemma 3.2, Proposition 3.3]{BL:Random_gen}. 

Denote by $E_a, E_b, E_{a,b}$ the closed fields generated by $\tr(a), \tr(b)$ and $(\tr(a),\tr(b))$ respectively. Where $\tr$ is the trace of the adjoint representation. Now for any given closed proper subfield $F < K$ the set $\{a \in A \ | \ \tr(a) \in F \}$ is of measure zero as demonstrated in \cite[Corollary 3.6]{BL:Random_gen}. Since there are only countably many finite sub fields in $K$ it follows that $E_a$ a local subfield of $K$, and hence $[K:E_a] < \infty$, for all but a nullset of elements. Therefore for almost all $a$ there are but finitely many intermediate subfields $E_a \le F \lneqq K$ and for every such $F$ the set $\{b \in A \ | \ \tr(b) \in F \}$ has Haar measure zero. The statement of the theorem now follows from Fubini's theorem. 
\end{proof}
\begin{remark}
Note that as stated the above criterion does not hold over fields of characteristic $2$. In fact it is easy to verify that even for the whole group $\tr(\PSL_2(K)) \subset K^2$ when $K$ is a field of characteristic $2$. Barnea and Larsen overcome this problem by passing to the simply connected cover and replacing the adjoint representation with the standard $2$-dimensional representation of $\SL_2$. 
\end{remark}

\subsection{Shalom's criterion} 
Using Pink's theorem, Shalom generalizes the criterion of openness to subgroups that are not compact in $\PSL_2$. We bring here a version of Shalom's theorem that is especially suited for our needs. As the proof is not difficult we repeat it for the convenience of the reader. 
\begin{proposition} (Shalom \cite[Proposition 8.4]{Shalom:Rigidity_comm_products}) \label{prop:Shalom_criterion}
Let $G$ be a non-discrete, unbounded closed subgroup of $\PSL_2(K)$. Assume that $G$ is Zariski dense and that $\tr(G) < K$ generates $K$ as a closed field. Then $G = \PSL_2(K)$.
\end{proposition}
\noindent
Before proceeding to the proof we recall the following
\begin{definition}
Let $G$ be a group and $H < G$ a subgroup, the {\it{commensurator}} of $H$ in $G$ is the subgroup $$\Comm_G(H) := \left\{g \in G \ | \ H \cap H^g {\text{ is of finite index in both }} H {\text{ and }} H^g \right\}.$$
Where $H^g = g H g^{-1}$.  We will also say that elements of $\Comm_G(H)$ {\it{commensurate the subgroup}} $H$.
\end{definition}
\noindent
Clearly $\Comm_G(H) = G$ for every compact open subgroup of $G$. 
Å\begin{proof}
It is well known that the only unbounded open subgroup of $\PSL_2(K)$ is the group itself (see for example \cite[Theorem T]{Prasad:Theorem_BTR}). So it is enough to prove that $G$ is open. Since $G$ is by assumption non-discrete and closed the group $M := G \cap \PSL_2(\Oc)$ is an infinite compact subgroup. Now by Pink's criterion it would be enough to show that $M$ is Zariski dense and that $\tr(M)$ is not contained in any closed proper subfield of $K$.

We first argue that $M$ is Zariski dense. Let $L := \left(\overline{M}^{Z}\right)^{0}$ be the connected component of the Zariski closure of $M$.  The group $L$ does not depend on finite index changes in $M$ and therefore it is normalized by $\Comm_{\PSL_2(K)}(M)$. In particular it is normalized by the group $G < \Comm_{\PSL_2(K)}(M) < N_{\PSL_2(K)}(L)$. Since $G$ is Zariski dense we conclude that $L$ is a normal subgroup of $\PSL_2$ and by simplicity of the latter group $L = \PSL_2$.  

Since $\chr(K) \ne 2$ the adjoint representation of $\PSL_2(K)$ is irreducible and we may assume use the simple form of Pink's Theorem (\ref{thm:Pink}): there exists a closed subfield $E < K$ and a model $H$ of $\PGL_2$ over $E$, such that $M$ is an open subgroup $H(E)$.  We argue that $G < H(E)$ using, again, the fact that $G$ commensurates $M$. Indeed for any $g \in G$, conjugation by $g$ takes a Zariski dense subgroup of $H(E)$ to a Zariski dense subgroup of $H(E)$ and therefore conjugation by $g$, and hence also $\Ad(g)$, are defined over $E$, for every $g \in G$. But $H$ is an adjoint group so the homomorphism $\Ad:G \arrow Ad(G)$ is an isomorphism, consequently $g \in H(E)$ for every $g \in G$. But hence $E = F$ and $\tr(G) \subset \tr(H(E))$ is contained in proper closed subfield of $K$, contradicting our assumption on the traces. 

We conclude that $M$ is open in $H(K) = \PGL_2(K)$ and consequently $G = \PSL_2(K)$. This completes the proof of the theorem. 
\end{proof}

\subsection{A subgroup of $\PSL_2(K)$ generated by an elliptic and a hyperbolic is generically dense}
\label{sec:SL_gen}
Collecting all of the above information we obtain the following analogue of Theorem \ref{thm:AG}
\begin{theorem} \label{thm:SL_gen} 
Let $K$ be a non-Archimedean local field of characteristic different than $2$ and $A = \PSL_2(K)$. Consider the set
$$\Hyp \times \Ell = \{(a,b) \in A^2 \ |  \ a {\text{ is hyperbolic and }} b {\text{ elliptic}} \}.$$  
Then $\Hyp \times \Ell$ is of infinite Haar measure and almost every $(a,b) \in \Hyp \times \Ell$ generates a dense subgroup of $A$.  
\end{theorem}
\begin{proof}
Both sets of hyperbolic and elliptic are open and invariant under conjugation and hence of infinite Haar measure. Let $G = \overline{\langle a, b \rangle}$ be the closed group generated by $a,b$. Because $a$ is, by assumption, hyperbolic the group $G$ cannot be bounded. Furthermore by \cite[Lemma 3.2]{BL:Random_gen} the Zariski closure of the cyclic subgroup $\langle b \rangle$ is almost surely a maximal torus. In particular the group $\overline{\langle b \rangle}$ is almost surely infinite (and compact) and hence $G$ cannot be discrete. The result now follows directly from Propositions \ref{prop:BL_conditions} and \ref{prop:Shalom_criterion}.
\end{proof}

\section{Proof of the main theorem}
From here on the proof for the group $A = \PSL_2(K)$ and for $A = \Aut^0(T)$ is identical, the only difference is that we will appeal to Theorem \ref{thm:SL_gen} for the former group and to Theorem \ref{thm:AG} for the latter. In particular by these two theorems the group $A$ is  topologically generated by two elements. We therefore consider the set $\Hom(F_k,A)$ for $k \ge 3$, with the action of $\Gamma = \Aut(F_k)$ by pre-composition. We will demonstrate that the $\Gamma$ action is ergodic on the invariant set $$D = D_{A,k}  \defeq \{\eta \in \Hom(F_k,A) \ | \ \overline{\eta(F_k)} = A\}$$ of homomorphisms with a dense image. For this proof it will be convenient to identify $\Hom(F_k,A)$ with $A^k$. To do this for every subset $I \subset \{1,2,\ldots,k\}$ we consider the map:
\begin{eqnarray*}
 \Psi_I: \Hom(F_k,A) & \arrow & A^I  \\
       \eta & \mapsto & \{\eta(x_i)\}_{i \in I},
\end{eqnarray*}
the map $\Psi_{\{1,2,\ldots,k\}}: \Hom(F_k,A) \arrow A^k$ is the natural isomorphism between these two groups. If we identify $\Hom(F_k,A)$ with $A^k$ the maps $\Psi_I: \Hom(F_k,A) \arrow A^I$ are just the projections on the corresponding partial products. It is easy to visualize the action of $\Gamma$ on $A^k$ defined by this identification. For example if 
$\a = (a_1,a_2,\ldots a_k) = \Psi_{\{1,2,\ldots,k\}}(\eta) = (\eta(x_1),\eta(x_2), \ldots, \eta(x_k))$, the Nielsen transformation $R_{i,j}$ will act by $R_{i,j} (a_1, \ldots, a_k) = (a_1, \ldots, a_i, \ldots, a_j a_i ,\ldots, a_k)$. In order to simplify the notation we will often identify $\Hom(F_k,A)$ with $A^k$ failing to mention the isomorphism $\Psi$ explicitly. 

Assume by way of contradiction that there is an essentially invariant subset $C \subset D$ which is neither null nor conull. Since $\Gamma$ is countable we may assume without loss of generality that $C$ is really invariant - replacing it if necessary by $\cap_{\gamma \in \Gamma} C^{\gamma}$. For $1 \le i \ne j \le k$ consider the set 
\begin{eqnarray*}
O_{i,j} &  \defeq & \{ \eta \in \Hom(F_k,A) \ | \ \eta(x_i) {\text{ is elliptic, and }} \eta(x_j) {\text{ is hyperbolic}} \}\\
& = & \{\a \in A^k \ | \ a_i \in \Ell, a_j \in \Hyp\}.
\end{eqnarray*}
For simplicity of notation we study the properties of $O \defeq O_{1,2} = \Ell \times \Hyp \times A^{k-2}$, but by symmetry all other sets $O_{i,j}$ will share the following nice properties:
\begin{enumerate}
\item \label{itm:open} $O$ is open of infinite measure in $\Hom(F_k,A)$.
\item \label{itm:inD} $O$ is essentially contained in $D$, i.e. $\mu(O \setminus D) = 0$.
\item \label{itm:FD} $O$ contains a fundamental domain for the action of $\Gamma$ on $D$. 
\end{enumerate}
For (\ref{itm:open}) use the fact that the sets of elliptic and hyperbolic elements are open sets of infinite Haar measure inside $A$. Theorem \ref{thm:AG}, or \ref{thm:SL_gen} as the case may be, directly implies (\ref{itm:inD}). Finally to demonstrate (\ref{itm:FD}) we have to find for every $\eta \in D$ an element $\gamma \in \Gamma$ such that $\eta \compos \gamma \in O$. The group $\eta(F_k)$ is dense in $A$ and in particular it does not act freely on $T$. Therefore, by Theorem \ref{thm:Weidmann}, for every $\eta \in D$ we can find some $\delta \in \Gamma$ such that $\eta \compos \delta(x_1)$ is elliptic. If $\eta \compos \delta (x_i)$ is hyperbolic for some $2 \le i \le k$ the element $\eta \compos \delta \compos T_{2,i} \in O$. If $\eta \compos \delta (x_i)$ is elliptic for all $i$ then, using the fact that $\eta \compos \delta(F_k)$ is dense and in particular does not fix a point on the tree, by \cite[p90, Corollaire 2]{Serre:Arbres} there exists some $2 \le i \le k$ such that $\eta \compos \delta(x_i x_1)$ is hyperbolic. Therefore $\eta \compos \delta \compos T_{i,2} \compos R^{+}_{1,i} \in O$. 

Since $O$ contains a fundamental domain and $C$ is $\Gamma$ invariant we can write:
\begin{equation} \label{eqn:FD}
C = \cup_{\gamma \in \Gamma} (C \cap O)^{\gamma}.
\end{equation}
It follows that $C \cap O$ inherits the non-triviality of $C$: it is neither null nor conull in $O$. 

We define 
$$E_{1,2} = \{(a,b) \in \Ell \times \Hyp \ | \ \mu^{k-2}(\Psi_{\{1,2\}}^{-1}((a,b)) \cap C) \ne 0\},$$  the set of all pairs $(a,b) \in \Ell \times \Hyp$ the fiber over which intersects $C$ in a non null subset of $A^{k-2}$. By Fubini's theorem $E_{1,2}$ is not null as a subset of $A^2$. 

\begin{figure}[htbp]
\begin{center}
\includegraphics[width=14cm]{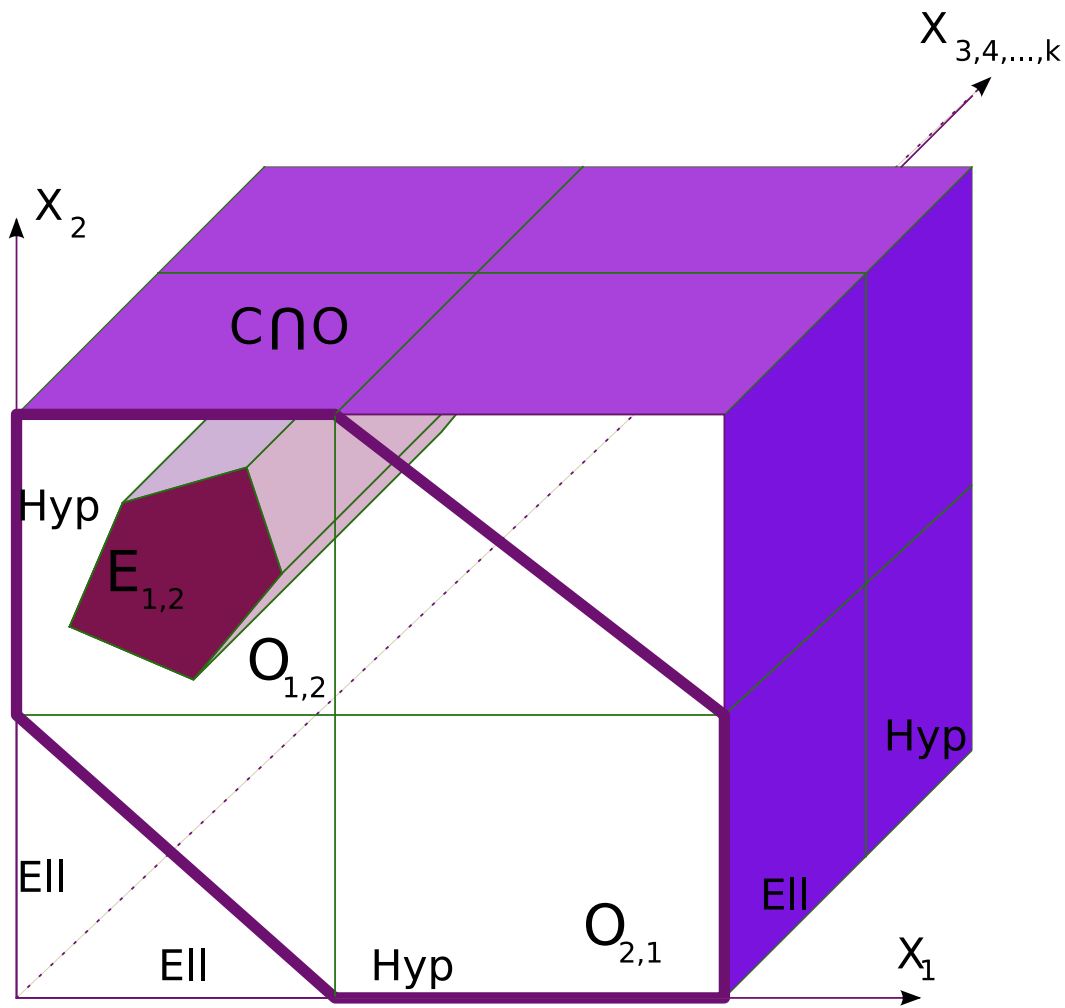}
\caption{The hexagon bounded by a thick line, on the front face, signifies the pairs of points that generate a dense subgroup. By Theorems \ref{thm:AG} or \ref{thm:SL_gen} this hexagon contains $O_{1,2}$ and $O_{2,1}$ - the upper left and lower right squares respectively. Over $O=O_{1,2}$ the set $C \cap O$ looks like a product.}
\label{fig:set_C}
\end{center}
\end{figure}

We now argue that for almost all $(a,b) \in E_{1,2}$ the fiber $C \cap \Psi_{1,2}^{-1}((a,b))$ is actually conull in $A^{k-2}$; by showing that these fibers are invariant under an ergodic group action. Indeed, the set $O \cap C$ is invariant under the subgroup of $\Gamma$ generated by the Nielsen transformations
 $\Sigma = \langle R^{\pm}_{1,i},R^{\pm}_{2,i} \ | \ 3 \le i \le k \rangle,$ 
and this action preserves the fibers over $\Psi_{1,2}$. Identifying the fiber $\Psi_{1,2}^{-1}((a,b))$ with a $A^{k-2}$, this action restricted to the fiber is just as the action of the subgroup $\langle a, b \rangle^{k-2}$ on $A^{k-2}$ by right multiplication. In particular whenever $\langle a,b \rangle < A$ is a dense subgroup this action is ergodic, by \cite[Lemma 2.2.13]{Zimmer:book}. This happens for almost every $(a,b) \in E_{1,2} \subset \Ell \times \Hyp$ using either Theorem \ref{thm:AG} or Theorem \ref{thm:SL_gen}. So up to nullsets we have identified the subset $O_{1,2} \cap C$ as 
\begin{equation} \label{eqn:product}
O_{1,2} \cap C = E_{1,2} \times A^{k-2}.
\end{equation}
In other words, the set $O_{1,2} \cap C$ is, measurably, a product. By the symmetry of the whole situation Equation (\ref{eqn:product}) is correct for any permutation of the index set. This situation is depicted in figure \ref{fig:set_C}. 

In the case $k \ge 4$ we can now conclude the proof. Equation (\ref{eqn:product}) for the coordinates $3,4$ says that $C \supset A^2 \times E_{3,4} \times A^{k-4}$. In particular for almost every $(a,b) \in A^2$ we have $\Psi_{1,2}^{-1}((a,b)) \supset E_{3,4} \times A^{k-4}$. But, just like $E_{1,2}$ we know that $\mu^2 (E_{3,4}) \ne 0$ so that $\mu^{k-2} (\Psi_{1,2}^{-1}((a,b))) \cap C) \ne 0$. In particular this is true for almost every $(a,b) \in \Ell \times \Hyp$. By definition, this means that $E_{1,2}$ is conull in $\Ell \times \Hyp$. Equation (\ref{eqn:product}) for the coordinates $1,2$ now reads $O_{1,2} \cap C = O_{1,2}$ so $O_{1,2}$ is essentially contained in  $C$. Since $C$ is $\Gamma$ invariant and $O_{1,2}$ contains a fundamental domain for the action of $\Gamma$ on $D$ we deduce that $C$ is conull in $D$ contradicting our definition of $C$. 

We now turn to the case $k = 3$. By definition of the product $\sigma$- algebra the set $E_{1,2}$ contains a cylinder $E_{1,2}^1 \times E_{1,2}^2$ where $E_{1,2}^1 \subset \Ell$ and $E_{1,2}^2 \subset \Hyp$ are non nullsets. Equation (\ref{eqn:product}) for coordinates $3,2$ (note that the order of the coordinates here does matter) reads $O_{3,2} \cap C = A \times E_{3,2}$. In particular for almost every $a \in A$ the set $\psi_{1}^{-1}(\{a\}) \cap O_{3,2} \cap C$ is essentially equal to $E_{3,2}$. Equation (\ref{eqn:product}) for coordinates $1,2$ shows that for a generic $a \in E_{1,2}^{1}$
$$E_{3,2} = \psi_{1}^{-1}(a) \supset E_{1,2}^2 \times \Ell.$$
 So that $\psi_{1}^{-1}(\{a\}) \cap O_{3,2} \cap C$ essentially contains $E_{1,2}^2 \times \Ell$ for almost every $a \in A$. Using Fubini we deduce that for almost every $(a,b) \in \Ell \times E_{1,2}^2$ the fiber $\psi_{1,2}^{-1}((a,b))$ is non null. So by the definition of $E_{1,2}$ we have $\Ell \times E_{1,2}^2 \subset E_{1,2}$. Repeating this argument with the indices $1,2$ exchanged and taking $E_{1,2}^1 = \Ell$ we conclude that $E_{1,2} = \Ell \times \Hyp$. So $O$ is essentially contained in $C$. But  $C$ is $\Gamma$ invariant and $O$ contains a fundamental domain for the $\Gamma$ action on $D$ so $C$ must be conull in $D$ which is a contradiction to our choice of $C$. This completes the proof of our main theorem.
\bibliography{../tex_utils/yair}
\end{document}